\input amstex
\input epsf
\documentstyle{amsppt}
\magnification=1095

\advance \hsize 2pt

\NoBlackBoxes

\define\RR{\Bbb{R}}

\define\ZZ{\Bbb{Z}}
\define\calL{{\Cal{L}}}

\def\eps{\varepsilon}
\def\ph{\varphi}
\def\del{\partial}

\let\ti=\tilde

\newcount\newrefnocount
\newrefnocount=1
\def\NewRefNo#1{\edef#1{\the\newrefnocount}\advance\newrefnocount1}

\NewRefNo\refAlex    
\NewRefNo\refBe      
\NewRefNo\refBi      
\NewRefNo\refBiMe    
\NewRefNo\refBiWr    
\NewRefNo\refChek    
\NewRefNo\refEli     
\NewRefNo\refElFr    
\NewRefNo\refEFM     
\NewRefNo\refFuTa    
\NewRefNo\refMarkov  
\NewRefNo\refMen     
\NewRefNo\refMor     
\NewRefNo\refWri     

\topmatter
\title          Markov theorem for transversal links 
\endtitle
\author         S.~Yu.~Orevkov, V.~V.~Shevchishin
\endauthor
\address 
S.~Yu.~Orevkov,
Steklov Math. Institute, ul. Gubkina 8, Moscow, Russia
\newline
Laboratoire E.Picard, UFR MIG, Universit\'e Paul Sabatier,
118 route de Narbonne, 31062 Toulouse, France
\endaddress
\email orevkov\@picard.ups-tlse.fr
\endemail
\address
V.~V.~Shevchishin,
Ruhr-Universit\"at, Bochum, 
Fakult\"at f\"ur Mathematik, 
Universit\"atsstrasse 150,
44801 Bochum, 
Germany
\endaddress
\email sewa\@cplx.ruhr-uni-bochum.de
\endemail
\abstract It is shown that two braids represent transversally isotopic links if
and only if one can pass from one braid to another by conjugations in braid
groups, {\it positive} Markov moves, and their inverses.

\bigskip
\centerline{ \hbox{\tenrm Revised version, 12 February 2002}}

\endabstract
\endtopmatter

\vskip-15pt

\document
By a well-known theorem of Alexander \cite{\refAlex}, any oriented link in
$\RR^3$ is isotopic to the closure of a braid. The question when two braids
represent isotopic links is answered by Markov's theorem \cite{\refMarkov} (see
\cite{\refBi}, \cite{\refBe}, or \cite{\refMor} for proofs): It is so if and
only if one can pass from one braid to another by conjugations in braid groups
$B_n$, the transformations $M_n^\pm: B_n \to B_{n+1}$, $M^{+}_n: b \mapsto
b\cdot \sigma_n$, $M^{-}_n: b \mapsto b\cdot \sigma^{-1}_n$ called {\it
positive/negative Markov moves} or {\it stabilizations}, and their inverses
({\it destabilizations}).

In the seminal paper \cite{\refBe} Bennequin established, among other very
important results, the analogue of Alexander's theorem for transversal
links (i.e., links transverse to the standard contact structure; see below). 
Namely, any transversal link is transversally isotopic to the closure of
a braid. The purpose of this paper is to prove the corresponding analogue of 
Markov's theorem.

\proclaim{Theorem} Two braids represent transversally isotopic links if and
only if one can pass from one braid to another by conjugations in braid
groups, {\it positive} Markov moves, and their inverses. 
\endproclaim

When this paper had been already finished, we learned from Victor Ginzburg
that he had announced this result around 1992. However, his proof has never
been published. Another proof of the theorem based on completely different
ideas was independently obtained by Nancy Wrinkle in her PhD thesis
\cite{\refWri}.

\medskip
Let us recall the standard definitions (see e.g.\ \cite{\refBe}). Consider the
1-form $\alpha= dz + x\,dy - y\,dx$ in $\RR^3$ with coordinates $x,y,z$. It
defines the standard contact structure in $\RR^3$. In the cylindric
coordinates $r,\theta, z$ with $x=r \cos \theta$, $y= r \sin \theta$ one has
$\alpha = dz + r^2 d\theta$.

A link $L$ in $\RR^3$ is {\it transversal} if the restriction $\alpha|_L$
nowhere vanishes. In this case $\alpha|_L$ defines a canonical orientation on
$L$.

A {\it geometric braid} in $\RR^3$ is an oriented link $L$ such that the
restriction $d\theta|_L$ is positive. In particular, $L$ is disjoint from the
$z$-axis $Oz$. The {\it degree} of $L$, also called the number of strings of
$L$, is the degree of the projection $(r, \theta, z) \mapsto \theta$
restricted to $L$. There is a canonical one-to-one correspondence between
isotopy classes of geometric braids of degree $n$ and conjugacy classes in the
braid group $B_n$.

Any conjugacy class in $B_n$ defines a transversal isotopy class of transversal
links. Indeed, any braid $b \in B_n$ can be realized as a geometric braid
sufficiently $C^1$-close to the standard circle $r=1$, $z=0$, which is clearly
transversal.

\smallskip 
The rest of the paper is devoted to the proof of Theorem. Essentially, our
proof is a parametric version of Bennequin's proof of his result cited above.

\smallskip 
Let $L_0$ and $L_1$ be two transversal geometric braids and $\{L_t \}_{t \in
[0,1]}$ a transversal isotopy between $L_0$ and $L_1$. Denote the interval
$[0,1]$ by $I$, the number of components of $L_0$ by $m$, and the disjoint
union of $m$ abstract circles by $S$. Abusing notation, we shall denote by
$s$\/ a positively oriented (local) coordinate on $S$, as also a current point
of $S$. The isotopy $\{L_t \}_{t\in I}$ can be parameterized by a smooth map
$\calL: S \times I \to \RR^3$ such that for every $t\in I$ the map $\calL_t: s
\mapsto \calL(s,t)$ is a parameterization of $L_t$.

\definition{Definition 1} Let $\{L_t \}_{t\in I}$ be a transversal isotopy
parameterized by a map $\calL: S \times I\to \RR^3$. It is called {\it
monotone near the axis} if there exists a finite number of parameters $0<t_1<
\dots < t_k <1$ such that the following holds:

\roster
\item For any $t_i$ there exists a unique $s_i \in S$ such that $\calL(s_i, t_i)$
lies on the $z$-axis $Oz$, and $\calL^{-1}(Oz)=\{(s_1,t_1),\dots,(s_k,t_k)\}$.
\item Up to a rotation of $\RR^3$ around $Oz$,
the mapping $\calL$ is given in a 
neighborhood of every $(s_i, t_i)$ by
$x = \tau - 3s^2$, $y = s\tau - s^3$, $z = z_i + s$, 
where $s$ is a positively oriented coordinate on $S$ centered at $s_i$ and
$\tau$ is a coordinate on $I$ centered at $t_i$ and oriented either positively
or negatively.
\endroster

The isotopy $\{L_t \}_{t\in I}$ is {\it monotone everywhere} if additionally
\roster
\item"(3)" 
$L_t$ is a transversal geometric braid for every $t \not\in
\{t_1,\ldots,t_k\}$.
\endroster

\enddefinition

Note, that if we fix $t\ne0$ and substitute $x = \tau - 3s^2$, 
$y = s\tau - s^3$ into
the $1$-form $r^2 d\theta=x\,dy-y\,dx$, we get 
$r^2d\theta =(\tau^2 +3s^4) ds>0$. 
Thus, conditions (2) and (3) of Definition 1 are consistent.

\smallskip
We shall always assume that isotopies we consider are sufficiently generic
outside a small neighborhood of the axis $Oz$.

\proclaim{ Lemma 1 } Let $b_0$ and $b_1$ be two braids, $L_0$ and $L_1$ the
transversal geometric braids defined by them. Assume that there exists an
everywhere monotone isotopy between $L_0$ and $L_1$. Then one can pass from
$b_0$ to $b_1$ by conjugations in braid groups, {\it positive} Markov moves,
and their inverses.
\endproclaim

\demo{Proof} When passing through a critical value $t=t_i$, the projection of
$L_t$ onto the horizontal plane $Oxy$ transforms near the origin as in Figure
1. This is a positive Markov move.
\qed
\enddemo

\midinsert 
\epsfxsize 70mm 
\centerline{\epsfbox{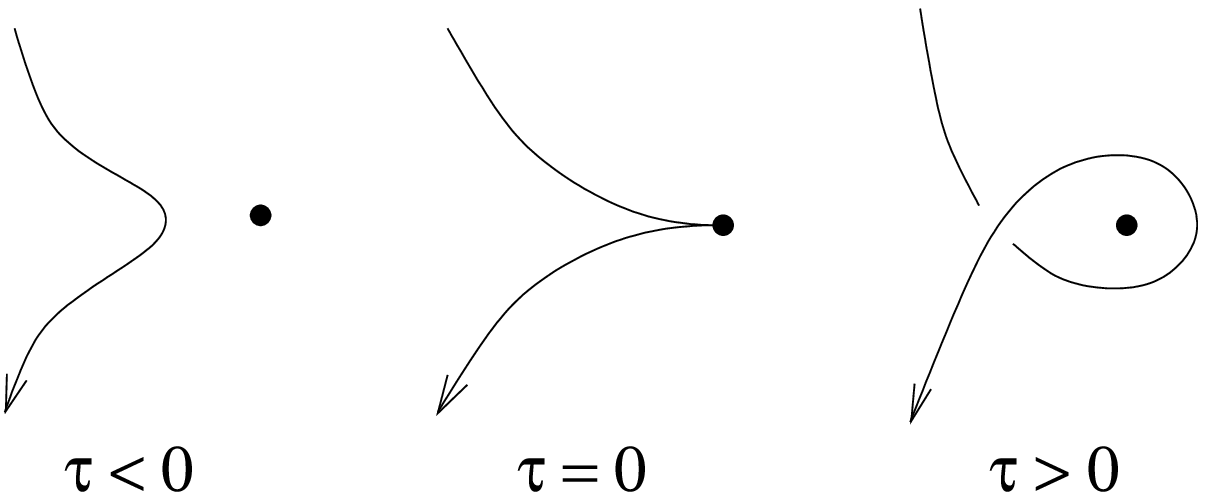}} 
\botcaption{ Figure 1. The curve $s\mapsto(\tau - 3 s^2, s\tau- s^3)$ } 
\endcaption 
\endinsert 

\bigskip

\proclaim{ Lemma 2 } Let $\{L_t \}_{t\in I}$ be a transversal isotopy
between transversal geometric braids $L_0$ and $L_1$. Then it can be perturbed
into an isotopy $\{L'_t \}_{t\in I}$ which is monotone near the axis.
\endproclaim

\demo{Proof} Let $\calL:S\times I\to\RR^3$ be a smooth mapping which
parameterizes $\{L_t\}$. Perturbing it if necessary, we can suppose that it is
transverse to the axis $Oz$. Let us consider a point $p=(s_0,t_0)\in S\times
I$ such that $\calL(p)$ lies on $Oz$. Let $s$ and $t$ be coordinates on $S$
and $I$ near $s_0$ and $t_0$ respectively (with $ds>0$). Set $\calL(s,t)=
\big(x(s,t), y(s,t), z(s,t)\big)$.  Since all $L_t$'s are transversal braids, 
we have $\del z/\del s>0$ at $p$. Hence, there exists a neighborhood $U$ of
$p$ such that $\del z/\del s>\varepsilon>0$ in $U$. Let us modify
$\big(x(s,t), y(s,t)\big)$ in $U$ replacing it by the homotopy in Figure 2
(the shaded zone corresponds to the homotopy described in Part (2) of
Definition 1 and shown in Figure 1; we assume here that before the
modification the homotopy looked as a parallel motion of a vertical line). If
$U$ is sufficiently small, then we can achieve that $|r^2\theta'_s| <
\varepsilon$ in $U$, which provides that $\calL^*_t\alpha>0$.
\qed
\enddemo

\midinsert
\epsfxsize 125mm
\centerline{\epsfbox{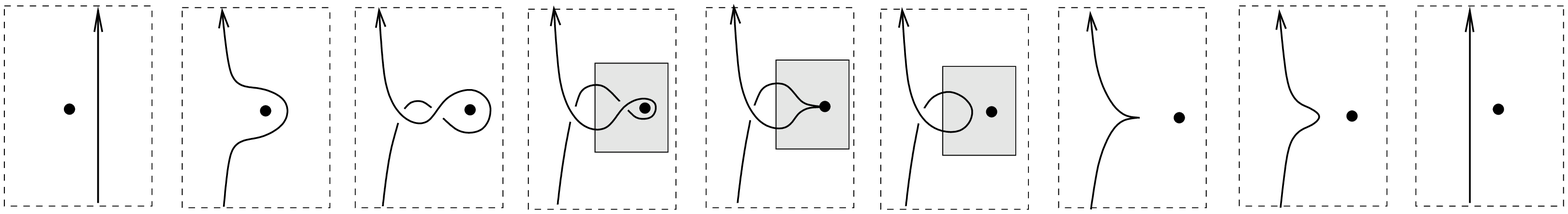}}
\botcaption{ Figure 2. Making the isotopy monotone near $Oz$ }
\endcaption
\endinsert

\smallskip

\definition{Definition 2} Let $\{L_t \}_{t\in I}$ be a transversal isotopy
parameterized by a map $\calL: S \times I\to \RR^3$. A {\it bad zone} of
$\calL$ is a connected component of the set of those points of $S \times I$ in
which $\del\theta/\del s \leqslant 0$, where $\theta(s,t)$ is the
$\theta$-component of $\calL(s,t)$.

\smallskip
A bad zone $V$ is {\it simple} if 
\roster
\item $V_t := (S \times t) \cap V$ is connected for all $t\in I$;
\item the total increment of $\theta$ along $V_t$ is less than $2\pi$.
\endroster

\smallskip 
The {\it shadow of $\calL$ on a bad zone $V$} is the set of those points $(s_0,
t_0) \in V$ for which the shortest segment connecting $p_0 := \calL(s_0,t_0)$
with the axis $Oz$ meets $L_{t_0}$ at some point $\calL(s_1,t_0)$.
The set of all such ``shading'' points $(s_1,t_0)$ will be called the 
{\it inverse shadow of $V$}.

\smallskip 
A bad zone $V$ is called {\it non-shadowed} if the shadow of $\calL$ on $V$
is empty.
\enddefinition


\medskip

\proclaim{ Lemma 3 } Let $\{L_t \}_{t\in I}$ be a transversal isotopy between
transversal geometric braids $L_0$ and $L_1$ parameterized by $\calL: S \times
I\to \RR^3$ which is monotone near the axis. Let $V$ be a simple and
non-shadowed bad zone and $U$ an arbitrary open subset of $S \times I$
containing $V$.

Then $\calL$ can be deformed into a transversal isotopy $\ti\calL: S \times
I\to \RR^3$ which is monotone near the axis, coincides with $\calL$ outside
$U$, and such that no bad zone of $\ti\calL$ meets $V$.
\endproclaim

\demo{Proof} Let us write in the cylindric coordinates $\calL(s,t) =
\big(r(s,t), \theta(s,t), z(s,t)\big)$. Then we have $z'_s + r^2\,\theta'_s
>0$. This implies that $z'_s >0$ on $V$. Choose a neighborhood $V^+$ of $V$
contained in $U$ such that $z'_s \geqslant \eps>0$ in $V^+$.

Let $[a,b]$ be the projection of $V$ onto $I$. We replace the components
$x(s,t)$ and $y(s,t)$ of $\calL$ in $V^+$ by the homotopy shown in Figure 3,
preserving the component $z(s,t)$.

\midinsert
\epsfxsize 115mm
\centerline{\epsfbox{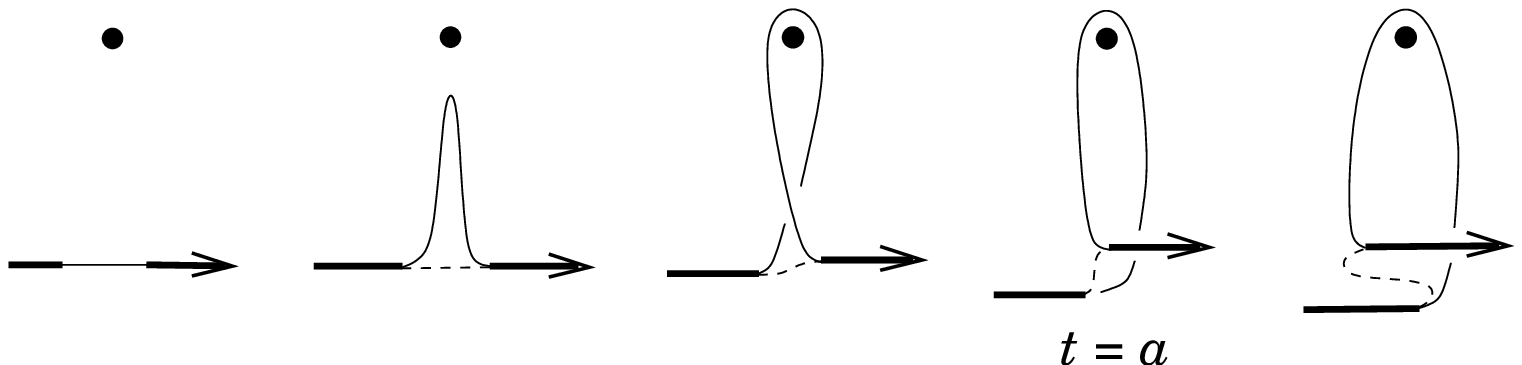}}
\botcaption{ Figure 3. Elimination of a bad zone (projection onto $Oxy$)}
\endcaption
\endinsert

In Figure 3, the bold lines represent the part of the homotopy which is not
changed; the dashed and resp.\ thin solid lines depict the isotopy before and
after the modification; the ``$\bullet$'' represents the origin of the plane
$Oxy$. The first three steps in Figure 3 is a deformation of the homotopy
described in Definition 1(2), see Figure 1.

Figure 3 depicts the modified homotopy for $t<c$ for some $c\in[a,b]$. To
construct the modified homotopy for $t>c$ we perform the same operations in the
reverse order.  \qed
\enddemo

\proclaim{ Lemma 4 } Let $\{L_t \}_{t\in I}$ be a transversal isotopy between
transversal geometric braids $L_0$ and $L_1$ parameterized by $\calL: S \times
I\to \RR^3$, which is monotone near the axis. Let $\big( r(s,t), \theta(s,t),
z(s,t) \big)$ be a representation of $\calL$ in cylindric coordinate. Let
$V$ be a bad zone, $l$ a generic smooth embedded curve in $V$ which is
the graph of a function $t=\ph(s)$,
and $U$ a
neighborhood of $l$ in $S \times I$. Let $\eps>0$.

Then there exist a sufficiently small open tubular neighborhood $U^-$ of $l$
in $S\times I$ and a perturbation $\ti\calL$ of $\calL$ of the form $\ti\calL
= \big( r(s,t), \ti\theta(s,t), z(s,t) \big)$ (i.e., only the
$\theta$-component is changed), such that

\roster
\item"(1)" $\del (V \backslash U^-)$ is smooth.
\item"(2)" $\ti\calL$ is monotone near the axis and coincides with $\calL$ 
outside $U$.
\item"(3)" $\del \ti\theta /\del s$ is positive in $U^-\cap V$ for $\ti\calL$.
\item"(4)" the signs of $\del \theta /\del s$ and $\del \ti\theta /\del s$ 
coincide outside $U^-\cap V$.
\item"(5)" $\max\limits_{U^-}\left({\del\ti\theta\over\del s}
\big/{\del z\over\del s}\right)
<\eps$. 

\endroster
\endproclaim

Informally speaking, this means that a bad zone can be cut along any smooth
curve. The operation described in the proof of Lemma 4 will be called {\it
wrinkling along the curve $l$}. The left hand side of (5) will be called the
{\it maximal slope of the wrinkling}. The assertion of the lemma in the
manifestation of the Gromov's $h$-principle in this setting.

\demo{ Proof }
In a neighborhood of every 
point $(s_0, t_0)$ of $l$ we perturb $\theta(s,t)$ by making a 
small wrinkle on the
graph of $\theta(s, t_0)$ at $s_0$ as it is shown in Figure 4, cf.\
\cite{\refBe}, pp.143--144.
\qed
\enddemo

\midinsert 
\epsfxsize 80mm 
\centerline{\epsfbox{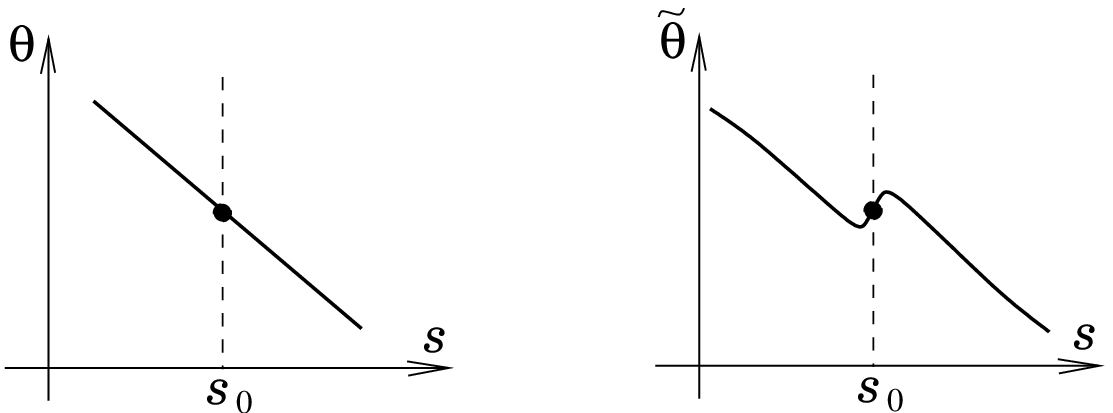}} 
\botcaption{ Figure 4. Wrinkling } 
\endcaption 
\endinsert 

\smallskip 
Let $\{L_t\}_{t\in I}$ be a transversal isotopy. Assume that $\{L_t\}$ is
monotone near the $Oz$-axis and generic outside a small neighborhood of the
axis $Oz$. Then for a generic value $t_0$ of the parameter $t$ the projection
of the link $L_{t_0}$ on the cylinder $S^1 \times \RR$ with the coordinates
$(\theta, z)$ is an immersion and the only singularities of the image are
crossings, i.e., ordinary double points. Moreover, there exist only finitely
many values $0<t_1 <\cdots< t_k<1$ for which the projection of $L_{t_i}$ on
$\theta z$-cylinder has a unique singularity of one of the following types:
\roster 
\item"(I)" $L_{t_i}$ meets the axis $Oz$ at some point in the way described
in Definition 1. 
\item"(II)" The projection of $L_{t_i}$ on $\theta z$-cylinder has a unique
ordinary tangency point.
\item"(III)" The projection of $L_{t_i}$ on $\theta z$-cylinder has a unique
ordinary triple point. 
\endroster
The singularities of types (II) and (III), respectively, are the second and
third Reidemeister moves in coordinates $(\theta, z)$. The first Reidemeister
move in coordinates $(\theta, z)$ is impossible for {\it transversal\/} links
since the derivatives ${\del z \over \del s}$ and ${\del \theta \over \del s}$
can not both vanish. 
Instead, a single Reidemeister move of the first 
kind occurs in every type (I) singularity of a transversal isotopy provided we 
consider the {\it projection on $Oxy$-plane}, see Figure 1. 

When we depict a crossing of the $\theta z$-projection of a link $L_t$, 
we assume that we look from the axis $Oz$, i.e. the overpass (resp. underpass) 
corresponds to the arc with a smaller (resp. bigger) value of $r$.
So, we say that an arc with a smaller value of $r$ {\it passes over} or
{\it shadows} an arc with a bigger value of $r$ (compare with Definition 2).

A singularity of the type (II) or (III) is called {\it positive} if 
${\del \theta \over \del s}>0$ 
at every point of $L_{t_i}$ which projects on the 
singularity, and {\it non-positive} otherwise.  
A non-positive singularity of the type (II) is called {\it bad} if 
there is a negative arc (with ${\del \theta \over \del s}>0$)
which is shadowed by another arc at the singularity.

\proclaim{ Lemma 5 }
Let $L$ be a transversal link. Suppose that the projection onto
the $\theta z$-cylinder has a bad non-positive 
singularity of the type {\rm(II)}.
Then the both branches are negative at this point.
\endproclaim

\demo{ Proof }
Let the branches be parametrized by 
$(r_\nu(s),\theta_\nu(s),z_\nu(s))$, $\nu=1,2$, so that $r_1>r_2$
at the crossing point. The tangency means that 
$z_2'/z_1'=\theta_2'/\theta_1'=\lambda$.
Since $\alpha|_L$ is positive, we have
$z_j' + r_j^2\theta_j'>0$, $j=1,2$.
Since the singularity is bad, we have $\theta_1'<0$. 
Suppose that $\theta_2'>0$. Then $\lambda<0$ and we have
$$
	0 < z_2' + r_2^2 \theta_2' < z_2' + r_1^2 \theta_2' =
	( z_1' + r_1^2\theta_1')\lambda < 0.\qquad\qed
$$
\enddemo

\proclaim{ Lemma 6 } Any transversal isotopy $\{L_t\}$ monotone near the
$Oz$-axis and generic outside it can be perturbed into a transversal isotopy
$\ti\calL$ without non-positive singularities of type {\rm(III)} and
without bad non-positive singularities of the type {\rm(II)}. 
 Moreover, such a perturbation can be made $C^0$-small and located
 in arbitrarily small neighborhoods of the points $(s_j, t_j)$ for which the
 thread $\calL(s, t_j)$ passes though a singularity of the type {\rm(II)} or
 {\rm(III)} with non-positive derivative ${\del \theta \over \del s}$ at
 $s=s_j$.
\endproclaim

\demo{ Proof } 
As in Lemma 4, it is sufficient to perturb only the coordinate $\theta$.
 
{\it Step 1. Elimination of non-positive triple points.} 
At each non-positive triple point, we perturb all negative branches as in 
Figure 5a. This can be done by replacing $\theta(s,t)$ with
$\tilde\theta(s,t)=\theta(s,t) + f(z(s,t),s)$ where the function $f(z,s)$
is the same for all the negative branches.
In the case when there are exactly two negative branches, we 
take care that for any $t$ the crossing point of the perturbed
branches rests on the same place as it was before the perturbation.
After such modification the triple point becomes positive and no other
triple points apear (a priori, new singularities of the type (II)
may appear).

{\it Step 2. Elimination of bad tangencies}.
Consider a bad non-positive singularity of the type (II).
By Lemma 5, the both branches are negative at this point.
We perturb them in the same way as in Step 1 (see Figure 5b).\qed

\midinsert 
\epsfxsize 120mm 
\centerline{\epsfbox{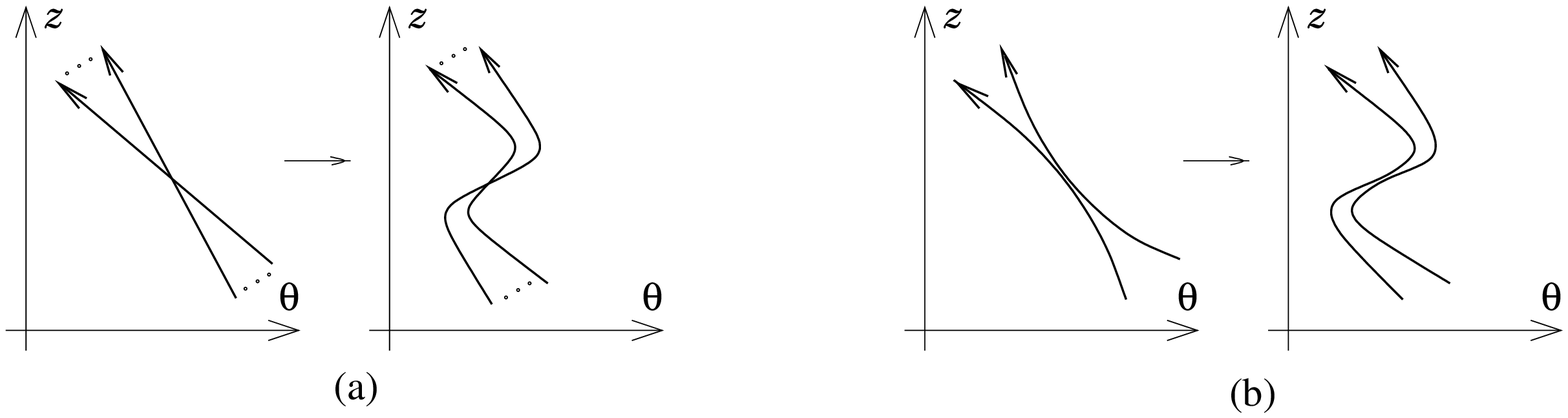}} 
\botcaption{ Figure 5. Elimination of bad non-positive singularities} 
\endcaption 
\endinsert 
 
\enddemo

\demo{Proof of Theorem} By Lemma 1, it is sufficient to prove that any
transversal isotopy $\calL$ between transversal geometric braids $L_0$ and
$L_1$ can be transformed into an everywhere monotone isotopy (see Definition
1). By Lemma 2, we may suppose that $\calL$ is monotone near the axis $Oz$.

Wrinkling $\calL$ along sufficiently many segments $s= \hbox{\it const}$ as in
Lemma 4, we can assume that all the bad zones are simple. Let us denote them
by $V_1, V_2, \ldots, V_n$. Fix disjoint neighborhoods $U_i$'s of $V_i$'s. We
are going to eliminate the bad zones one by one modifying $\calL$ at the
$i$-th step only in $U_i\cup\dots\cup U_n$. This insures that the procedure
will terminate. The isotopy obtained after the $i$-th modification is
denoted by $\calL_i$ and $\calL_0 = \calL$ is the initial isotopy. 
Every $\calL_i$
will be monotone near the axis $Oz$.

To pass from $\calL_i$ to $\calL_{i+1}$, we proceed as follows (compare with
\cite{\refBe}, Theorem 8, pp.142--144).

\roster
\item"{\it Step 1}." Eliminate non-positive singularities of
$\calL_i$ of the type (III) and bad non-positive singularities
of the type (II) applying Lemma 6.
\endroster

\noindent
Let us consider connected components $\ell_1,\ell_2,\dots$ 
of the inverse shadow of $V_i$
on the other bad zones (a bad zone cannot shadow itself because
$\partial z/\partial s>0$ on it). 
Any point $(s,t)$ of any $\ell_\nu$ corresponds to a 
crossing of the projection of $L_t$
onto the $\theta z$-cylinder. 
The crossing is either as in Figure 6a or as in Figure 6b.
 
\midinsert 
\epsfxsize 120mm 
\centerline{\epsfbox{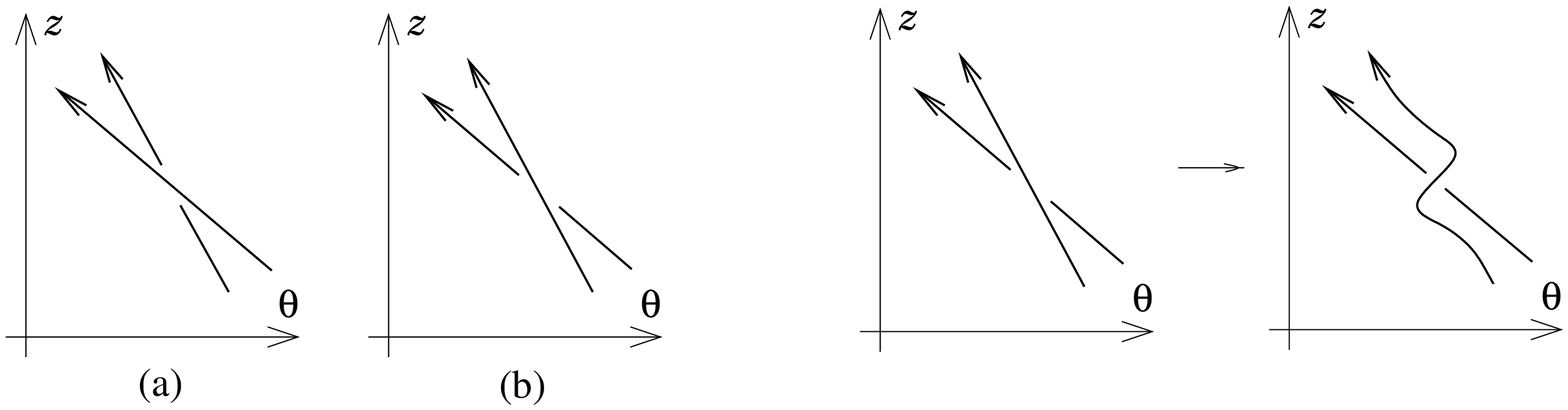}} 
\botcaption{ \hskip 10mm
		Figure 6.
	\hskip 25mm 
		Figure 7. Wrinkling at Step 2 } 
\endcaption 
\endinsert 

\roster
\item"{\it Step 2}."
For each component $\ell_\nu$ corresponding to Figure 6b, 
we wrinkle the corresponding bad zone along it (see Figure 7).
\item"{\it Step 3}." Wrinkle $V_i$ along the shadow of $\calL_i$
(see Figure 8).
\endroster

\midinsert 
\epsfxsize 120mm 
\centerline{\epsfbox{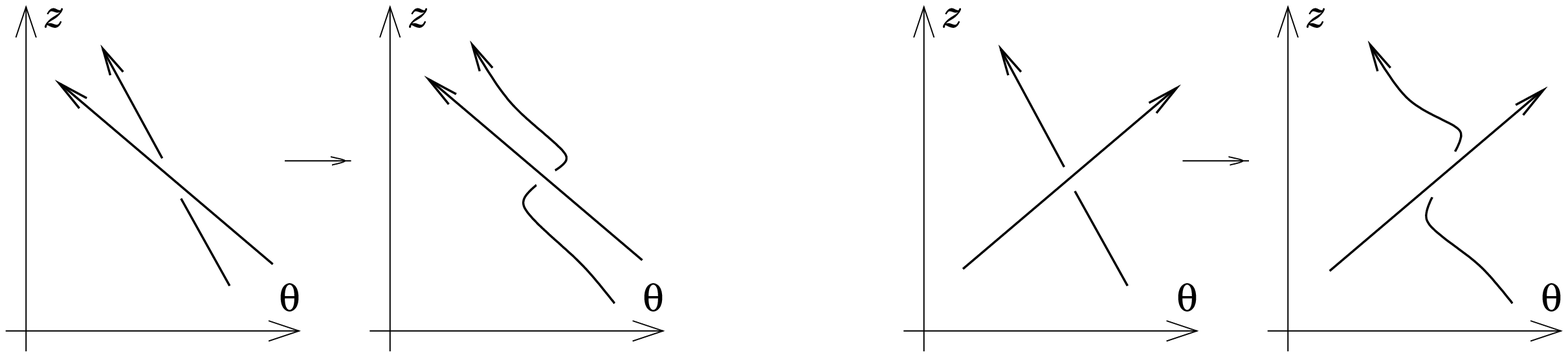}} 
\botcaption{ Figure 8. Wrinkling at Step 3 } 
\endcaption 
\endinsert 

\noindent
Note that crossings as in Figure 6a are eliminated at Step 2 and
the fact that crossings as in Figure 9 are impossible, is proved
in  \cite{\refBe, pp.142--144} (the proof is similar to that of Lemma 5).
If the maximal slope of the wrinkling is small enough (see condition (5) of
Lemma 4), then no new shadow appears because the wrinkling is performed
away from tangencies and triple points.

\midinsert
\epsfxsize 25mm
\centerline{\epsfbox{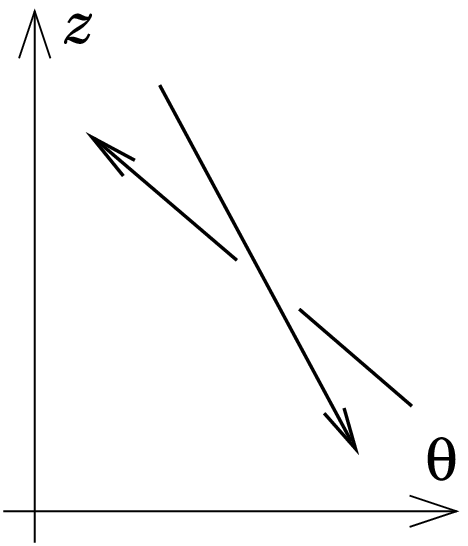}}
\botcaption{ Figure 9. Impossible crossing }
\endcaption
\endinsert

\roster
\item"{\it Step 4}." Wrinkle, if necessary, the obtained bad zones along
segments $s= \hbox{\it const}$ to make all the bad zones simple.
\item"{\it Step 5}." Apply Lemma 3 to all the newly obtained bad zones in $U_i$.
\qed
\endroster
\enddemo

\example{ Example } According to \cite{\refEli}, two transversal unknots are
transversally isotopic iff they have the same Bennequin index. The Bennequin
index of a transversal geometric braid $L$ corresponding to a braid $b \in
B_n$ is equal to $(\sum_i k_i)-n$ for $b = \prod_i \sigma_{j_i} ^{k_i}$ (see
\cite{\refBe}). Therefore, by our Theorem, any braid representing an unknot
can be transformed by positive (de)stabilizations and conjugations into the
braid $\sigma_1^{-1} \dots \sigma_{n-1}^{-1} \in B_n$ for some $n$. Here is
the sequence of transformations for the braid $\sigma_1^{-1} \sigma_2
\sigma_3^{-1}$ ($k$ and $\bar k$ stand for $\sigma_k$ and $\sigma_k^{-1}$;
$M_+^{-1}$ for a positive destabilization):
$$
\bar12\bar3  =  \bar1\bar332\bar3  = \bar3\bar132\bar3  = \bar3 \bar1 \bar2 3 2
\buildrel conj \over \longrightarrow 
\bar1\bar232\bar3  = \bar1\bar2\bar232  
\buildrel conj \over \longrightarrow 
 2\bar1\bar2\bar23  
\buildrel M_+^{-1} \over \longrightarrow 2\bar1\bar2\bar2 
\buildrel conj \over \longrightarrow 
 \bar1\bar2.
$$
\endexample

\head Appendix. Markov's Theorem from the Point View of Contact Topology.
\endhead

Here we discuss some ``classical'' and recent results on contact
isotopy of Legendrian and transversal knots in $\RR^3$ and deduce
``topological'' Markov's theorem from its contact version.

\smallskip
We start with a brief description of related notions and constructions,
referring to the articles \cite{\refEFM} and \cite{\refFuTa} for more details. 
Notice that the contact structure in $\RR^3$ use there is given by the form
$\alpha_{jet} := dz -ydx$ and originates from the identification of $\RR^3$
with the space $J^1\RR$ of $1$-jets of functions on the real axis $\RR$. The
substitution $(x,y,z) \mapsto (x, 2y, z+xy)$ transforms $\alpha_{jet}$ into
the rotation invariant form $\alpha_{rot} := dz +xdy - ydx= dz + r^2 d\theta$
used in the main part. Thus both forms define the same contact structure. The
advantage of the form $\alpha_{jet}$ is that it provides a possibility to
control over Legendrian and transversal knots by their projections on $xy$-
and $xz$-planes.

A link $L$ in $\RR^3$ is {\it Legendrian} ({\it transversal}) w.r.t.\ a
contact form $\alpha$ if the restriction $\alpha|_L$ vanishes identically
(never vanishes). The {\it contact} orientation of a transversal link $L$ is
induced by the restriction $\alpha|_L$. A link isotopy $\{L_t\}$ is Legendrian
(resp., transversal) if every $L_t$ has this property. We always assume that 
Legendrian, transversal, and usual (``topological'') isotopies preserve the
orientation of the link.

Every link in $\RR^3$ admits both Legendrian and transversal representation. 
Legendrian and transversal links have additional $\ZZ$-valued invariants
constraining the existence of a contact isotopy: these are the Maslov and
Thurston-Bennequin indices in the Legendrian case and the Thurston-Bennequin
index in the transversal case, denoted by $\mu(L)$ and $tb(L)$ respectively.

Assume that $L_1$ and $L_2$ are disjoint links, both Legendrian or
transversal. Let $lk(L_1, L_2)$ be their linking number. Then $\mu(L_1 \sqcup
L_2) = \mu(L_1) + \mu(L_2)$ (linear behavior) and $tb(L_1 \sqcup L_2) =
tb(L_1) + 2lk(L_1, L_2) + tb(L_2)$ (quadratic behavior). This reduces the
computation of the indices to the case of knots.

The Thurston-Bennequin index of a {\it knot} is independent of its
orientation, while the Maslov index changes the sign if we reverse the
orientation. The Thurston-Bennequin index of a transversal link $L(b)$
represented by an algebraic braid $b$ with $n$ strands equals $tb(L(b)) =
\deg(b)-n$ where $\deg(b)$ is the algebraic degree of $b$.

Every oriented Legendrian link $L$ can be smoothly approximated by a
transversal link $L^+$ whose contact orientation coincides with that induced
from $L$. Moreover, such a link $L^+$ is unique up to transversal isotopy. 
Similarly, there exists a unique transversal isotopy class of links $L^-$
which approximate $L$ with the reversed orientation. The indices of $L^\pm$
are related to those of $L$ as $tb(L^\pm) = tb(L_0) \pm \mu(L_0)$.

There exist several constrains on possible values of Maslov and
Thurston-Benne\-quin indices of Legendrian and transversal links in $\RR^3$. 
The first one is that $tb(L)$ (resp., $tb(L) \pm \mu(L)$) has the same parity
as the number of components of the transversal (Legendrian) link $L$. In
particular, $tb(L)$ is odd for every transversal knot. Another constrain is
the Bennequin inequality $tb(L) \leq -\chi(F)$ for every transversal link $L$
and its Seifert surface $F$. Unlike the first constrain, this one is highly
non-trivial and reflects the fact that the standard contact structure in
$\RR^3$ is {\it tight} (see \cite{\refEli} for more details). For a Legendrian
link $L$ this inequality reads $tb(L) + |\mu(L)| \leq -\chi(F)$. Some further
inequalities are listed in \cite{\refFuTa}.

It is always possible to decrease the Thurston-Bennequin index of a Legendrian
or transversal knot $L$. More precisely, there exists transformations
$\zeta_+$ and $\zeta_-$ (resp., a transformation $\rho$) of isotopy classes of
oriented Legendrian (resp., transversal) knots with the following properties:

\roster
\item"(1)" The transformations $\zeta_\pm$ and $\rho$ can be realized by
adding an appropriate unknotted loop in any given neighborhood of any given
point on $L$; in particular, they can be represented by an appropriate
Legendrian knot $\zeta_\pm L$ (resp., a transversal knot $\rho L$) in the
topological isotopy class of $L$.
\item"(2)" The operations $\zeta_+$ and $\zeta_-$ commute, i.e., there exists
a Legendrian isotopy between $\zeta_+ (\zeta_- L)$ and $\zeta_- (\zeta_+ L)$.
\item"(3)" If a braid $b$ represents a transversal knot $L$, then the braid
$M^-b$ obtained from $b$ by negative stabilization represents $\rho L$.
\item"(4)" $tb(\zeta_\pm L) = tb(L) -1$ and $\mu(\zeta_\pm L) = \mu(L) \pm 1$
in the Legendrian case;  $tb(\rho L) = tb(L) -2$ in the transversal case.
\item"(5)" For an oriented Legendrian knot $L$, the knot $(\zeta_+ L)^+$ (see
above) is transversally isotopic to $L^+$ and the knot $(\zeta_+ L)^-$ to
$\rho(L^-)$; similarly, the knot $(\zeta_- L)^-$ is transversally isotopic to
$L^-$ and the knot $(\zeta_- L)^+$ to $\rho(L^-)$.
\item"(6)" The transformations $\zeta_\pm$ and $\rho$ naturally extend to
links; one should only indicate to which component of the link the operation
is applied.
\endroster

We refer to \cite{\refFuTa} for the definition of the transformations
$\zeta_\pm$ and the proof of the properties (1--5). However, it should be
noticed that these transformations are known well enough as a part of the
contact topology folklore, so looking for references would be an ungrateful
task. The property (3) can be used as the definition of the operation $\rho$. 
The property (5) means that, informally speaking, after the ``positive
(negative) transversalization'' $L \mapsto L^+$ (resp., $L \mapsto L^-$) the
operation $\zeta_\pm$ descends to the stabilization $M^\pm$ (resp., $M^\mp$)
of the same (resp., opposite) sign.

\proclaim{ Proposition A}
Let $L_1$ and $L_2$ be two oriented Legendrian (resp., transversal)
links which are topologically isotopic; then one can obtain Legendrian (resp.,
transversal) isotopic links $L'_1$ and $L'_2$ applying the operations
$\zeta_\pm$ (resp., $\rho$) to each component of $L_1$ and $L_2$ sufficiently
many times.
\endproclaim

\proclaim{ Proposition B} 
Let $L_1$ and $L_2$ be two oriented Legendrian links; then the links $L_1^+$
and $L_2^+$ are transversally isotopic if and only if one can transform $L_1$
into $L_2$ applying Legendrian isotopies, the operation $\zeta_+$, and its
inverse.
\endproclaim

In the case of knots Proposition A was proved in \cite{\refFuTa} and
Proposition B in \cite{\refEFM}. However, since the condition of being
connected is not used in the both proofs, the general case follows as well. In
view of the property (3), our Theorem and Proposition A imply Markov's
theorem for knots in the refined form stated in Introduction. As one can
easily see the refined form remains valid in the case of links after an
appropriate generalization of negative (de)stabilizations. Such a
generalization should represent the operation $\rho$ applied to any prescribed
component of the link. For example, one can take operations 
$$ 
M^-_k: b \in B_n \;\mapsto \;
\sigma_{n-1}\cdots \sigma_k\, b\, \sigma_k^{-1}\cdots \sigma_{n-1}
^{-1} \sigma_n^{-1}  \in B_{n+1}
$$ 
which are compositions of the conjugation in $B_n$ by $\sigma_{n-1}\cdots
\sigma_k$ with the negative stabilization $M^-$.

\medskip
In view of Proposition A, the authors of \cite{\refFuTa} have expressed the
conjecture that the transversal (Legendrian) isotopy class of a knot is
completely determined by its topological isotopy class and its
Thurston-Bennequin (and Maslov) index. This conjecture has been disproved by
Yuri\u\i\ Chekanov who has constructed \cite{\refChek} new invariants of
Legendrian knots and has given an example of two Legendrian knots which are
topologically isotopic and have equal Thurston-Bennequin and Maslov indices
but different Chekanov's invariants. Some examples of even finer type have
been found in \cite{\refEFM}. Namely, there exist Legendrian knots $L_1$ and
$L_2$ which have equal Thurston-Bennequin and Maslov classes and
transversally isotopic ``transversalizations'' $L^+_1$ and $L^+_2$, but
nevertheless $L_1$ and $L_2$ are not Legendrian isotopic. On the other hand,
the Legendrian isotopy class of the unknot is completely determined by its
Thurston-Bennequin and Maslov indices, see \cite{\refEli} and \cite{\refElFr}.

A similar counterexample for transversal knots has been constructed in 
\cite{\refBiMe}. It is shown that the braids 
$$
\sigma_1^{2p+1} \sigma_2^{2q} \sigma_1^{2r} \sigma_2^{-1}
\quad\text{and}\quad
\sigma_1^{2p+1} \sigma_2^{-1} \sigma_1^{2r} \sigma_2^{2q}
\quad 
\text{ with }
p,q,r>1 \text{ and } q\ne r
$$
represent the knots $K_1$ and $K_2$ which are topologically isotopic and have
equal Thurston-Bennequin indices but which are not isotopic transversally.
On the other hand, there are several types of knots and links for which the
transversal isotopy class is completely determined by its topological isotopy
class and Thurston-Bennequin indices of the components, see  \cite{\refBiWr}.
For example, those are unlinks and iterated torus knots.

The discussions made so far lead to the following problems:
\roster
\item"{\bf P1}" {\sl Does there exist two Legendrian knots $L_1$ and 
$L_2$ which are not Legendrian isotopic, but the ``transversalizations'' of
both signs $L^+_1$ and $L^+_2$ (resp., $L^-_1$ and $L^-_2$) are transversally
isotopic?}
\endroster
The negative answer to this question is conjectured (indirectly) in
\cite{\refEFM}. 
\roster
\item"{\bf P2}" {\sl Find an analogue of Alexander's and Markov's theorems
for Legendrian links}. 
\endroster

\smallskip
We finish the paper with a description of a natural construction of closed
Legendrian braids. It can be considered as the first step toward the solution
of Problem P2.
First, we describe possible Legendrian isotopy classes of unknots. Let $\bar
L_{0,0}$ be the curve in the $xz$-plane given by the equation $z^2 =\cos^3(x)$
with $|x| \leq \pi/2$ and $|z| \leq1$, see Figure 10.
\midinsert
\epsfxsize 100mm
\centerline{\epsfbox{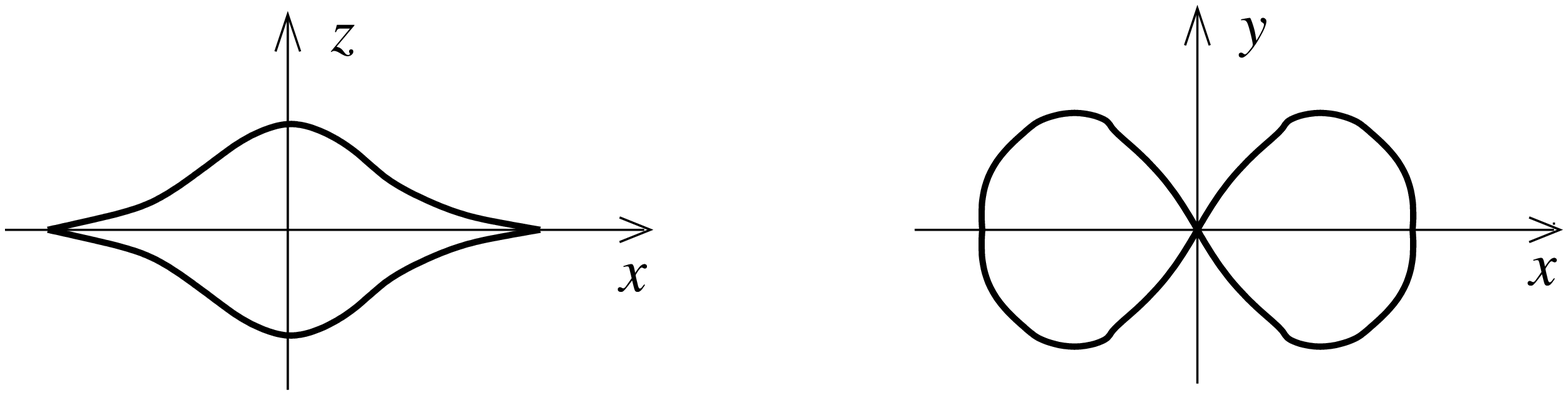}}
\botcaption{ Figure 10. The $xz$- and $xy$-projections of $L_{0,0}$.}
\endcaption
\endinsert
This curve lifts uniquely to a smooth Legendrian curve $L_{0,0}$ in $\RR^3$
with the standard contact structure given by $dz -ydx$. Namely, the lift of
each branch given $z(x)= \pm \cos^{3/2}(x)$ is parameterized by $(x, z'(x),
z(x))$ with $z'(x) := {dz(x) \over dx} = \mp {3\over 2} \cos^{1/2}(x)
\sin(x)$. Observe that in  a neighborhood of each cusp-point $(\pm{\pi\over2},
0, 0)$ the curve $L_{0,0}$ admits the parameterization 
$$
x(t) = \pm \arccos(t^2),\quad
y(t) = \mp {\textstyle{3\over 2}} t \sqrt{1-t^4},\quad
z(t) = t^3
$$ 
with $t$ close to $0$. This shows that $L_{0,0}$ is a smooth Legendrian
unknot. Direct computation gives $\mu(L_{0,0}) =0$ and $tb(L_{0,0}) = -1$,
see \cite{\refFuTa}, \S3.4. 

Set $L(p,q) := \zeta^p_+ \circ \zeta^q_- (L_{0,0})$. Then $\mu(L_{p,q}) =p-q$
and $tb(L_{p,q}) = -1-p-q$. By the results of Bennequin and
Eliashberg that every Legendrian unknot $L$ is Legendrian isotopic to $L(p,q)$
with $p = (\mu(L) - tb(L)-1)/2$ and $q = (-\mu(L) - tb(L)-1)/2$.

\medskip
Now assume that $L_b$ is a Legendrian braid with the ``axis'' $L_a$, which is
also a Legendrian knot. Then there exists a tubular neighborhood $U\cong
\Delta \times L_a$ of $L_a$ and coordinates $(v,w;\theta)$ in $U$ such that
\roster
\item $\theta$ is the coordinate along $L_a \cong S^1$;
\item $(v,w) \in \Delta$ are normal coordinates to $L_a$ in $U$;
\item the contact structure in $U$ is given by the form $dv -wd\theta$;
\item the projection $\pi_{v,\theta}: L_b \to [0,1] \times L_a$ of $L_b$ onto
the strip $[0,1] \times L_a$ with coordinates $(v,\theta)$ has only simple
transversal crossings. 
\endroster 
Observe that the projection of $L_b$ onto $(v,\theta)$-strip determines $L_b$
completely. Indeed, if $(v, w) =(f_i(\theta), g_i(\theta))$ is a local
parameterization of a strand of $L_b$, then $g_i(\theta)$ is the derivative of
$f_i(\theta)$, $g_i(\theta) = f'_i(\theta)$. It follows then that the
projection has only {\it positive} crossings. 

Vice versa, given a Legendrian knot $L_a$ and a positive braid $b$, there
exists a Legendrian link $L_b$ realized as the closure of $b$ in arbitrary
tubular neighborhood $U$ of $L_a$. Moreover, the Legendrian isotopy class of
such a link $L_b$ is well-defined. We shall use the notation $L_a \ltimes b$
to denote such a link $L_b$.

\proclaim{ Lemma 6} 
\roster
\item The link $L_{p,q} \ltimes b$ is represented by the braid
$\Delta^{-p-q-1} \cdot b$. 
\item For any Legendrian knot $L_a$ and a positive braid $b\in B_+(n)$, 
\item" "
$\mu\big(L_a\ltimes b \big) = \mu(L_a) \cdot n
\qquad\text{and}\qquad
tb\big(L_a \ltimes b\big) = n^2 tb(L_a) + \deg(b)
$.
\endroster
In particular, $\mu\big(L_{p,q}\ltimes b \big) = n(p-q)$ and
$ tb\big(L_{p,q} \ltimes b\big) = \deg(b) -n^2(p+q+1)$.
\endproclaim

Observe that every braid $b \in B(n)$ can be decomposed as $b= \Delta^{-k}
\cdot b_+$ with appropriate $k\ge0$ and $b_+\in B_+(n)$.

\demo{ Proof } Every Legendrian knot $L$ in $\RR^3$ has two natural framings:
the Legendrian one given by the contact distribution $\xi := ker(dz -ydx)
\subset T\RR^3$ and the topological one given by its Seifert surface. In
particular, the coordinates $(v,w, \theta)$ in a tubular neighborhood of $L$
introduced above define the Legendrian framing. By definition, the
Thurston-Bennequin index $tb(L)$ is the linking number between $L$ and the
knot $L'$ obtained from $L$ by pushing it slightly in the positive (or
negative) normal direction to the contact distribution $\xi = ker(dz -ydx)$. 
Thus $tb(L)$ is the rotation number of the Legendrian framing with respect to
the topological one. The part $(1)$ of the lemma follows.

It follows from definition that the Maslov index of a Legendrian link $L$
in $\RR^3$ is the winding number of the projection of $L$ onto $xy$-plane.
Since every strand of $L_a \ltimes b$ is $C^1$ close to $L_a$ we immediately
obtain $\mu\big(L_a\ltimes b \big) = \mu(L_a) \cdot n$.

Now assume that $b_0 \in B_+(n)$ is the trivial braid. Let $L_i$, $i=1\ldots
n$, be the strands of $L_a \ltimes b_0$. Then every $L_i$ is Legendrian
isotopic to $L_a$ and represents the ``push in the direction normal to
$\xi$''. So the linking number $lk(L_i, L_j) = tb(L_i) = tb(L_a)$. Then
$tb(L_a \ltimes b_0) =\sum_i tb(L_i) + \sum_{i<j} 2 lk(L_i, L_j) = n^2
tb(L_a)$. 

To obtain the general case, we use the algorithm for computing of the
Thurston-Bennequin index of a Legendrian link $L$ in $\RR^3$ by its projection
onto $xz$-plane, see \cite{\refFuTa}, \S3.4. After a small Legendrian
perturbation, the only singularities of such a projection are transversal
crossings and cusps. A crossing is called positive (negative) if both strands
cross the vertical line in the same (resp., opposite) direction. Then $tb(L)$
is the number of positive crossings minus the number of negative crossings
minus half the number of cusps. Now, it remains to observe that for $b\in
B_+(n)$ the $xz$-projections of $L_a \ltimes b$---compared with that of $L_a
\ltimes b_0$---has exactly $\deg(b)$ additional positive crossings.
\qed
\enddemo

\smallskip\noindent
{\bf Acknowledgement.} The authors were supported by Deutsche
Forschungsgemeinschaft Schwerpunkt ``Global Methods in Complex Geometry''.

\enddocument